\documentclass[12pt]{amsart}

\usepackage{amsmath,amsfonts,amscd,amsthm,amssymb,latexsym}
\usepackage[all,cmtip]{xy}

\font\teneufm=eufm10 \font\seveneufm=eufm7 \font\fiveeufm=eufm5
\newfam\frakturfam
\textfont\frakturfam=\teneufm \scriptfont\frakturfam=\seveneufm
\scriptscriptfont\frakturfam=\fiveeufm

\newtheorem{pr}{Proposition}

\newtheorem{lm}{Lemma}
\newtheorem{theor}{Theorem}
\newtheorem{co}{Corollary}

\def\bee{\begin{eqnarray}}
\def\bes{\begin{eqnarray*}}
\def\eee{\end{eqnarray}}
\def\ees{\end{eqnarray*}}

\def\a{\alpha}
\def\b{\beta}
\def\g{\gamma}

\def\Proof{{\sl Proof.}\ }

\title{Automorphisms of affine Veronese surfaces}

\begin{document}
\date{}
\maketitle

\begin{center}
{\bf Bakhyt Aitzhanova}
\footnote{Department of Mathematics,
 Wayne State University,
Detroit, MI 48202, USA, e-mail: {\em aitzhanova.bakhyt01@gmail.com}}
and
{\bf Ualbai Umirbaev}\footnote{Department of Mathematics,
 Wayne State University,
Detroit, MI 48202, USA 
and Institute of Mathematics and Mathematical Modeling, Almaty, 050010, Kazakhstan,
e-mail: {\em umirbaev@wayne.edu}}

\end{center}

\begin{abstract} We prove that every derivation and every locally nilpotent derivation of the subalgebra $K[x^n, x^{n-1}y,\ldots,xy^{n-1}, y^n]$, where $n\geq 2$, of the polynomial algebra $K[x,y]$ in two variables over a field $K$ of characteristic zero is induced by a derivation and a locally nilpotent derivation of $K[x,y]$, respectively. Moreover, we prove that every automorphism of  $K[x^n, x^{n-1}y,\ldots,xy^{n-1}, y^n]$ over an algebraically closed field $K$ of characteristic zero is induced by an automorphism of $K[x,y]$. We also show that the group of automorphisms of $K[x^n, x^{n-1}y,\ldots,xy^{n-1}, y^n]$ admits an amalgamated free product structure. 
\end{abstract}


\noindent


\section{Introduction}

\hspace*{\parindent}

Let $K$ be an arbitrary field and let $\mathbb{A}^n$ and $\mathbb{P}^n$ be the affine and the projective $n$-space over $K$, respectively. 
The \textit{Veronese map} of degree $d$ is the map
\bes
\nu_d:\mathbb{P}^n\to\mathbb{P}^m
\ees 
that sends $[x_0:\ldots:x_n]$ to all $m+1$ possible monomials of total degree $d$, where 
\bes
m={n+d \choose d}-1.
\ees
It is well known that the image of the Veronese map is a projective variety and is called the \textit{Veronese variety} \cite{Har}. 

The \textit{rational normal curve} $C_n\subset \mathbb{P}^n$ is a particular case of the Veronese variety and  is defined to be the image of the map
\bes
\nu_n: \mathbb{P}^1\to\mathbb{P}^n
\ees
given by
\bes
\nu_n:[x_0:x_1]\mapsto[x_0^n:x_0^{n-1}x_1:\ldots:x_1^n]=[X_0:\ldots:X_n].
\ees

It is well known that $C_n$ is the common zero locus of the polynomials 
\bee\label{f1}
F_{i,j}=X_{i-1}X_{j+1}-X_iX_j \,\, \textit{for} \,\,1\leq i\leq j\leq n-1.
\eee

For $n = 2$ it is the plane conic $X_0X_2=X_1^2$  and for $n= 3$ it is the twisted cubic  \cite{Har}.

Denote by $V_n\subset \mathbb{A}^{n+1}$  the common zero locus of the polynomials  (\ref{f1}) in  $\mathbb{A}^{n+1}$.   The variety $V_n$ is the affine cone of the rational normal curve $C_n$ and is called the Veronese cone in \cite{KPZ2017}. We will call $V_n$ the \textit{affine Veronese surface of index $n$}  in order to separate it from the Veronese cones of higher dimensions.  Veronese surfaces play an important role in the description of quasihomogeneous affine surfaces given by M.H. Gizatullin \cite{Giz1971} and V.L. Popov \cite{Popov1973}. They form one of the main examples of the so called {\em Gizatullin surfaces} \cite{KPZ2017}. 

 M.H. Gizatullin and V.I. Danilov  devoted two papers \cite{GD1975,GD1977}  to the systematic study of automorphisms of affine surfaces including affine cones of rational normal curves. In particular, generators of the automorphism group of $V_n$ can be deduced from their work along with its amalgamated product structure. L. Makar-Limanov \cite{ML90,ML01} gave an algebraic description of generators of the automorpism groups of algebraic surfaces defined by an equation of the form $x^ny=P(z)$. This gives an explicit description of generators  of the automorphism group of $V_2$.

It is well known \cite{Jung, Kulk} that all automorphisms of the polynomial algebra $K[x,y]$ in two variables $x,y$ over a field $K$ are tame. 
The well-known Nagata automorphism (see \cite{Nagata})
\bes
\sigma=(x+2y(zx-y^2)+z(zx-y^2)^2, y+z(zx-y^2),z)
\ees
of the polynomial algebra $K[x,y,z]$ over a field $K$ of characteristic zero is proven to be non-tame \cite{SU}.

The automorphism group $\mathrm{Aut}\,K[x,y]$ of this algebra admits an amalgamated free product structure \cite{Kulk, Shafar}, i.e.,
\bee\label{f0}
\mathrm{Aut}\,K[x,y]=\mathrm{Aff}_2\,(K)\ast_{C}\mathrm{Tr}_2\,(K),
\eee
where $\mathrm{Aff}_2\,(K)$ is the group of affine automorphisms,  $\mathrm{Tr}_2\,(K)$ is the group of triangular automorphisms, and  $C=\mathrm{Aff}_2\,(K)\cap \mathrm{Tr}_2\,(K)$.

It follows that any algebraic subgroup $G\subseteq \mathrm{Aut}\,K[x,y]$ is conjugate to a subgroup of one of the factors 
$\mathrm{Aff}_2\,(K)$ and $\mathrm{Tr}_2\,(K)$ \cite{GD1977,Kam1979,Wright1979}. In particular, any reductive subgroup $G\subseteq \mathrm{Aut}\,K[x,y]$ is {\em linearizable}, i.e., is conjugated to a subgroup of linear automorphisms $GL_2(K)$.  
The first examples of nonlinearizable actions  were given by G.W. Schwarz \cite{Schwarz89} and a nonlinearizable action of the symmetric group $S_3$ on $\mathbb{C}^4$ is given in \cite{FM02}. It is still open question if any finite automorphism of $\mathbb{C}^n$ for $n\geq 3$ is linearizable \cite{KS1995}.

Recently I. Arzhancev and M. Zaidenberg \cite{AZ2013} proved that every automorphism of the Veronese surface $V_n$ can be extended to an automorphism of the plane using the construction of Cox rings. It is also shown that the automorphism group of the Veronese surface $V_n$ admits an amalgamated product structure induced by (\ref{f0}) and an analogue of the Kambayashi \cite{Kam1979} and Wright \cite{Wright1979} result for $V_n$ is proven.

This paper is devoted to the study of vector fields and automorphisms of the affine Veronese surface $V_n$ for all $n\geq 2$ by purely algebraic methods. The algebra of polynomial functions on $V_n$ is isomorphic to the subalgebra 
$K[x^n, x^{n-1}y,\ldots,xy^{n-1}, y^n]$ of $K[x,y]$ (Proposition \ref{p0}). Thus the group of automorphisms of $V_n$ is anti-isomorphic to the group of automorphisms of the algebra 
$K[x^n, x^{n-1}y,\ldots,xy^{n-1}, y^n]$. We show that over a field $K$ of characteristic zero every derivation  and every locally nilpotent derivation of the algebra $K[x^n, x^{n-1}y,\ldots,xy^{n-1}, y^n]$ is induced by a derivation  and a locally nilpotent derivation of $K[x,y]$, respectively. Using the proof of the Rentchler's Theorem \cite{RR} on locally nilpotent derivations of $K[x,y]$ given in \cite[Ch. 5]{Es}, we prove that every automorphism of $K[x^n, x^{n-1}y,\ldots,xy^{n-1}, y^n]$ is induced by an automorphism of $K[x,y]$ if $K$ is an algebraically closed field of characteristic zero.  This gives an explicit description of generators of the automorphism group of $V_n$ as opposed to papers \cite{AZ2013,GD1977}. 
We also show that the amalgamated free product structure of the automorphism group of $K[x,y]$ induces an amalgamated free product structure on the automorphism group of $K[x^n, x^{n-1}y,\ldots,xy^{n-1}, y^n]$.

The paper is organized as follows. In Section 2 we describe the algebra of polynomial functions on the affine Veronese surface $V_n$. 
In Section 3 we recall some necessary results on the structure of the automorphism group of $K[x,y]$ from \cite{Cohn,Es}. Section 4 is devoted to lifting of derivations of $K[x^n, x^{n-1}y,\ldots,xy^{n-1}, y^n]$ to derivations of $K[x,y]$. In Section 5 we prove that so called $n$-graded derivations of $K[x,y]$ are triangulable. In Section 6 we prove that every automorphism of $K[x^n, x^{n-1}y,\ldots,xy^{n-1}, y^n]$ is induced by an automorphism of $K[x,y]$. The amalgamated free product structure of the automorphism group of $K[x^n, x^{n-1}y,\ldots,xy^{n-1}, y^n]$ is given in Section 7.

 \section{Polynomial functions on $V_n$}

\hspace*{\parindent}

Let $K$ be an arbitrary field and let $K[X_0,X_1,\ldots,X_n]$ be the polynomial algebra over $K$ in the variables $X_0,X_1,\ldots,X_n$.  The set of all monomials of the form
\bee \label{3}
u=X_0^{i_0}X_1^{i_1}\ldots X_n^{i_n},
\eee
where $ i_0,i_1,\ldots,i_n\geq0$, is a linear basis of $K[X_0,X_1,\ldots,X_n]$. Set $\a(u)=(i_0,i_1,\ldots,i_n)\in \mathbb{Z}^{n+1}$. 
If $u$ and $v$ are two monomials of the form (\ref{3}) then set $u\leq v$ if $\a(u)\leq \a(v)$ with respect to the lexicographical order. 

Let $I$ be the ideal of $K[X_0,X_1,\ldots,X_n]$ generated by all elements $F_{ij}$ defined in (\ref{f1}).

\begin{lm} \label{l0} The images of all different monomials of the form $X_k^iX_{k+1}^j$, where $0\leq k\leq n-1$ and $i,j\geq 0$, in $K[X_0,X_1,\ldots,X_n]/I$ form a linear basis of $K[X_0,X_1,\ldots,X_n]/I$. 
\end{lm}
\Proof The leading monomial of $F_{ij}$ with respect to the ordering $\leq$ is $X_{i-1}X_{j+1}$. Consider the leading monomials 
$X_{i-1}X_{j+1}$ and $X_{k-1}X_{l+1}$ of $F_{ij}$ and $F_{kl}$, respectively. Assume $i\leq k$ and $F_{ij}\neq F_{kl}$. Then monomials $X_{i-1}X_{j+1}$ and $X_{k-1}X_{l+1}$ have nontrivial intersection in the following cases: 

(a) $i=k$ and $j<l$; 

(b) $j+1=k-1$; 

(c) $i<k$ and $j=l$. 

Case (a). We form an $S$-polynomial (see, for example \cite{CLO})
\bes
S(F_{ij},F_{il})=(X_{i-1}X_{j+1}-X_iX_j)X_{l+1}-(X_{i-1}X_{l+1}-X_iX_l)X_{j+1}\\
=-(X_jX_{l+1}-X_{j+1}X_l)X_{i}=-F_{j+1,l}X_i.
\ees
The leading monomial of $F_{j+1,l}X_i$ is equal to $X_{i}X_jX_{l+1}$ and is less than $X_{i-1}X_{j+1}X_{l+1}$. 

Case (b). We have $i+2\leq j+2\leq l$. Then 
\bes
S(F_{ij},F_{(j+2)l})=(X_{i-1}X_{j+1}-X_iX_j)X_{l+1}-X_{i-1}(X_{j+1}X_{l+1}-X_{j+2}X_l)\\
=-X_iX_jX_{l+1}+X_{i-1}X_{j+2}X_l=F_{i(j+1)}X_l-X_iF_{(j+1)l}. 
\ees 
The leading term of $F_{i(j+1)}X_l$ is $X_{i-1}X_{j+2}X_l$ and the leading term of $X_iF_{(j+1)l}$ is $X_iX_jX_{l+1}$. Both of them are less than $X_{i-1}X_{j+1}X_{l+1}$. 

Case (c) can be handled similar to the case (a). 

Consequently, the set of all elements $F_{ij}$ forms a Gr\"{o}bner basis  for the ideal $I$ \cite[Theorem 6, p. 86]{CLO}. Since the leading monomial of $F_{ij}$  is $X_{i-1}X_{j+1}$ it follows the statement of the lemma \cite[Ch. 5, section 3]{CLO}. $\Box$

\begin{pr} \label{p0}
$K[X_0,X_1,\ldots,X_n]/I\cong K[x^n,x^{n-1}y,\ldots,xy^{n-1},y^n]$. 
\end{pr}
\Proof The homomorphism $$\phi:K[X_0,X_1,\ldots,X_n]\to K[x^n,x^{n-1}y,\ldots,xy^{n-1},y^n]$$ defined by 
$\phi(X_i)=x^{n-i}y^i$ for all $i$ induces the homomorphism 
$$\bar\phi: K[X_0,X_1,\ldots,X_n]/I\to K[x^n,x^{n-1}y,\ldots,xy^{n-1},y^n]$$ since $\phi(X_{i-1}X_{j+1}-X_iX_j)=x^{n-(i-1)}y^{i-1}x^{n-(j+1)}y^{j+1}-x^{n-i}y^ix^{n-j}y^j=0$  for all $1\leq i\leq j\leq n-1$. 

Let 
$u=X_k^iX_{k+1}^j$ and $v=X_s^pX_{s+1}^q$ where $k\leq s$. We get  
\bes
\phi(u)=(x^{n-k}y^k)^i(x^{n-k-1}y^{k+1})^j=x^{(n-k)i+(n-k-1)j}y^{ki+(k+1)j}
\ees
and, similarly, 
\bes
\phi(v)=x^{(n-s)p+(n-s-1)q}y^{sp+(s+1)q}. 
\ees
Consequently, $\phi(u)=\phi(v)$ if and only if 
\bee\label{b1}
\nonumber (n-k)i+(n-k-1)j=(n-s)p+(n-s-1)q, \\
 ki+(k+1)j=sp+(s+1)q. 
\eee
By adding the both sides of these equalities we get $n(i+j)=n(p+q)$, i.e., $i+j=p+q$. Then (\ref{b1}) gives that 
\bes
k(p+q)+j=s(p+q)+q, 
\ees
i.e., 
\bee\label{b2}
(s-k)(p+q)=j-q. 
\eee
We get $j-q\geq 0$ since $s\geq k$. Then (\ref{b2}) 
 is possible only if $s=k$ or $s-k=1$ and $p+q=j-q$.  If $s=k$ then (\ref{b2}) gives $j=q$. Then $i=p$ since  $i+j=p+q$. This gives $u=v$. Suppose that $s-k=1$ and $p+q=j-q$. Since $i+j=p+q$ it follows that $q=0, i=0,p=j$. Then $u=v=X_{k+1}^j$. 

  Thus we proved that the images of different monomials of the form $X_k^iX_{k+1}^j$ under $\phi$ are different monomials in $x,y$. Consequently, the images of different monomials of the form $X_k^iX_{k+1}^j$ are linearly independent. 

By Lemma \ref{l0}, the images of different monomials of the form $X_k^iX_{k+1}^j$ gives a linear basis for $K[X_0,X_1,\ldots,X_n]/I$. Consequently, $\bar\phi$ is an injection.  
Obviously, $\bar\phi$ is a surjection, i.e., $\bar\phi$ is an isomorphism. 
$\Box$

 \section{Automorphisms of $K[x,y]$}

\hspace*{\parindent}

Let $K[x,y]$ be the polynomial algebra in the variables $x,y$ over a field $K$ and let $\mathrm{Aut}\,K[x,y]$ be the group of automorphisms of  $K[x,y]$. Denote by $\phi=(f,g)$ the automorphism of $K[x,y]$ such that $\phi(x)=f$ and $\phi(y)=g$, where $f,g\in K[x,y]$. If $\phi=(f_1,g_1)$ and $\psi=(f_2,g_2)$ then the product in $\mathrm{Aut}\,K[x,y]$ is defined by
$$\phi\circ\psi=(f_2(f_1,g_1),g_2(f_1,g_1)).$$

An automorphism $\phi\in\mathrm{Aut}\,K[x,y]$ is called \textit{elementary} if it has the form
$$\phi=(x,\alpha y+f(x))$$
or
$$\phi=(\alpha x+g(y),y),$$
where $f(x)\in K[x]$, $g(y)\in K[y]$,  and $0\neq\alpha\in K$. The subgroup of $\mathrm{Aut}\,K[x,y]$ generated by all elementary automorphisms is called the \textit{tame subgroup}. Elements of this subgroup are called \textit{tame automorphisms} of $K[x,y]$. 

An automorphism $\phi\in\mathrm{Aut}\,K[x,y]$ is called \textit{affine} if it has the form
$$\phi=(\alpha_1x+\beta_1y+\gamma_1,\alpha_2x+\beta_2y+\gamma_2)$$
where $\alpha_1\beta_2\neq\beta_1\alpha_2$ and $\alpha_1,\alpha_2,\beta_1,\beta_2,\gamma_1,\gamma_2\in K$. The subgroup $\mathrm{Aff}_2\,(K)$ of $\mathrm{Aut}\,K[x,y]$ generated by all affine automorphisms is called the \textit{affine subgroup}.  If $\gamma_1,\gamma_2=0$ then the affine automorphism $\phi$ is called \textit{linear}. The subgroup $\mathrm{GL}_2\,(K)$ of $\mathrm{Aff}_2\,(K)$ generated by all linear automorphisms is called the \textit{linear subgroup}. 

An automorphism $\phi\in\mathrm{Aut}\,K[x,y]$ is called \textit{triangular} if it has the form
\bee\label{f2}
\phi=(\alpha x+f(y),\beta y + \gamma),
\eee
where $0\neq\alpha,\beta\in K$ and $f(y)\in K[y]$. The subgroup $\mathrm{Tr}_2\,(K)$ of $\mathrm{Aut}\,K[x,y]$ generated by all triangular automorphisms is called the  \textit{triangular subgroup}.

The well known Jung-van der Kulk Theorem \cite{Jung,Kulk} says that all automorphisms of the polynomial algebra $K[x,y]$ in two variables $x,y$ over a field $K$ are tame. Moreover,  van der Kulk and Shafarevich  \cite{Kulk, Shafar} proved that the automorphism group  $\mathrm{Aut}\,K[x,y]$ of this algebra admits an amalgamated free product structure, i.e.,
\bes
\mathrm{Aut}\,K[x,y]=\mathrm{Aff}_2\,(K)\ast_{C}\mathrm{Tr}_2\,(K),
\ees
where $C=\mathrm{Aff}_2\,(K)\cap \mathrm{Tr}_2\,(K)$.

We fix a grading
\bee\label{g}
K[x,y]=K[x,y]_0\oplus K[x,y]_1\oplus K[x,y]_2 \oplus \ldots\oplus K[x,y]_{n-1}
\eee
of the polynomial algebra $K[x,y]$, where $K[x,y]_i$ is the linear span of all homogeneous monomials of degree $i+ns$, $i=0,1,\ldots,n-1$, and $s$ is an arbitrary nonnegative integer. This is a $\mathbb{Z}_n$-grading of $K[x,y]$, i.e., 
$$K[x,y]_iK[x,y]_j\subseteq K[x,y]_{i+j},$$ where $i,j\in  \mathbb{Z}_n=\mathbb{Z}/n\mathbb{Z}$. For shortness we will refer to this grading as $n$-grading.

An automorphism $\phi\in\mathrm{Aut}\,K[x,y]$ is called a {\em graded automorphism} with respect to grading  (\ref{g}) if $\phi(x), \phi(y)\in K[x,y]_1$. A graded automorphism is called {\em graded tame} if it is a product of graded elementary automorphisms. 

Recently A. Trushin \cite{Trushin} studied graded automorphisms of polynomial automorphisms. But his gradings do not include gradings of type  (\ref{g}). 

A graded automorphism  of $K[x,y]$ with respect to grading  (\ref{g}) will be called an {\em n-graded automorphism} for shortness.  Obviosly, every $n$-graded automorphism induces an automorphism of the algebra 
 $K[x^n,x^{n-1}y, \ldots,xy^{n-1}, y^n]$. 

A  derivation $D$ of $K[x,y]$ will be called an {\em n-graded derivation} if $D(x), D(y)\in K[x,y]_1$. Recall that every derivation $D$ of $K[x,y]$ can be uniquely written in the form 
\bes
D=f\partial_x+g\partial_y, 
\ees
where $D(x)=f$, $D(y)=g$, and $\partial_x=\frac{\partial}{\partial x}$ and $\partial_y=\frac{\partial}{\partial y}$ are partial derivatives with respect to $x$ and $y$, respectively.

\section{Derivations of $K[x^n, x^{n-1}y,\ldots,xy^{n-1}, y^n]$}

\hspace*{\parindent}

Let $K$ be an arbitrary field of characteristic zero. 
Let $A$ be any algebra over $K$. A derivation $D$ of $A$ is  called {\em locally nilpotent} if for every $a\in A$ there exists a positive integer $n=n(a)$ such that $D^n(a)=0$. 

If $D$ is a locally nilpotent derivation of  $A$ then 
\bes
\exp{D}=\sum_{p\geq 0} \frac{1}{p!}D^p
\ees
is an automorpism of $A$ and  is called an {\em exponential} automorphism. 

Moreover, if $D$ is any derivation of  $A$ then 
\bes
\exp{TD}=\sum_{i=0}^{\infty} \frac{1}{i!}D^iT^i
\ees
is an automorpism of the formal power series algebra $A[[T]]$.  If $D$ is locally nilpotent then $\exp{TD}$ is an automorphism of $A[T]$.

Consider the grading (\ref{g})  of $K[x,y]$. 
A  derivation $D$ of $K[x,y]$ will be called an {\em n-graded derivation} if $D(x), D(y)\in K[x,y]_1$. 
Obviously, every $n$-graded derivation of $K[x,y]$ induces a derivation of $K[x^n, x^{n-1}y,\ldots,xy^{n-1}, y^n]$. The reverse is also true.

\begin{lm}\label{l1}
Every derivation of $K[x^n, x^{n-1}y,\ldots,xy^{n-1}, y^n]$ can be uniquely extended to an $n$-graded derivation of $K[x,y]$. 
\end{lm}
\Proof Let $D$ be a derivation of $K[x^n, x^{n-1}y,\ldots,xy^{n-1}, y^n]$. Denote by $T$ the unique extension of $D$  \cite[p. 120]{ZS} to a derivation of the field of fractions $K(x^n, x^{n-1}y,\ldots,xy^{n-1}, y^n)$ of 
$K[x^n, x^{n-1}y,\ldots,xy^{n-1}, y^n]$. Obviously, the field extension 
$$K(x^n, x^{n-1}y,\ldots,xy^{n-1}, y^n)\subseteq K(x,y)$$
 is algebraic. This extension is separable since $K$ is a field of characteristic zero. By Corollaries 2 and 2' in \cite[pages 124--125]{ZS}, every derivation of the field 
$K(x^n, x^{n-1}y,\ldots,xy^{n-1},$\\$ y^n)$ can be uniquely extended to a derivation of $K(x,y)$. Let $S$ be the unique extension of $T$ to a derivation of $K(x,y)$. Suppose that 
\bee\label{f3}
S(x)=\frac{f_1}{g_1}, \, S(y)=\frac{f_2}{g_2},
\eee
where $f_1,f_2\in K[x,y]$, $0\neq g_1,g_2\in K[x,y]$, and the pairs $f_1,g_1$ and $f_2,g_2$ are relatively prime. 
We have
\bes
D(x^{n-i}y^i)=S(x^{n-i}y^i)=(n-i)x^{n-i-1}y^i\frac{f_1}{g_1}+ix^{n-i}y^{i-1}\frac{f_2}{g_2}
\ees
for all $0\leq i\leq n$.

 Since $D(x^n), D(x^{n-1}y),\ldots,D(xy^{n-1}), D(y^n)\in K[x^n, x^{n-1}y,\ldots,xy^{n-1}, y^n]$ it follows that 
\bes
g_1g_2|(n-i)x^{n-(i+1)}y^if_1g_2+ix^{n-i}y^{i-1}f_2g_1
\ees
for all $0\leq i\leq n$. Consequently, 
\bes
g_1|(n-i) x^{n-(i+1)}y^{i}
\ees
and
\bes
g_2| i x^{n-i}y^{i-1}
\ees
for all $0\leq i\leq n$. 

This means that $g_1|x^{n-1}$ and $g_1|y^{n-1}$ and, consequently, we may assume that $g_1=1$.
Similarly, $g_2|y^{n-1}$ and $g_2|x^{n-1}$ give that $g_2=1$. Obviously, $f_1,f_2\in K[x,y]_1$. $\Box$\\

For any derivation $D$ of $K[x^n, x^{n-1}y,\ldots,xy^{n-1}, y^n]$ denote by $\widetilde{D}$ its unique extension to a derivation of $K[x,y]$ determined by Lemma \ref{l1}. Obviously $D$ is locally nilpotent if $\widetilde{D}$ is locally nilpotent. The reverse statement is also true.

\begin{lm}\label{l2}
If $D$ is a locally nilpotent derivation of $K[x^n, x^{n-1}y,\ldots,xy^{n-1}, y^n]$ then $\widetilde{D}$ is a locally nilpotent $n$-graded derivation of $K[x,y]$. 
\end{lm}
\Proof
Suppose that $D$ is a locally nilpotent derivation of $K[x^n, x^{n-1}y,\ldots,xy^{n-1}, y^n]$. Then $\exp{TD}$ is an automorphism of $K[x^n, x^{n-1}y,\ldots,xy^{n-1}, y^n][T]$. Recall that $\exp{T\widetilde{D}}$ is an automorphism of  $K[x,y][[T]]$.    We have 
\bes
\exp{TD}(x^n)=\exp{T\widetilde{D}}(x^n)=\exp{T\widetilde{D}}(x)^n. 
\ees
This implies that $\exp{T\widetilde{D}}(x)\in K[x,y][T]$ since $\exp{TD}(x^n)\in K[x,y][T]$. Similarly, $\exp{T\widetilde{D}}(y)\in K[x,y][T]$. This means that there exist natural numbers $m$ and $n$ such that $\widetilde{D}^m(x)=0$ and $\widetilde{D}^n(y)=0$. Therefore $\widetilde{D}$ is locally nilpotent. 
$\Box$\\

\section{Triangulation of locally nilpotent $n$-graded derivations}

\hspace*{\parindent}

A derivation $D$ of $K[x,y]$ is called {\em triangular} if 
\bes
D(x)=f(y)\in K[y], \ \ D(y)=\a \in K. 
\ees

A derivation  $D$ of $K[x,y]$ is called {\em triangulable} if there exists an automorphism $\alpha\in\mathrm{Aut}\,K[x,y]$ such that $\alpha^{-1}D\alpha$ is triangular. 

Every triangular derivation, and hence every triangulable derivation, is locally nilpotent.  In 1968 R. Rentschler \cite{RR} proved that every locally nilpotent derivation of the polynomial algebra $K[x,y]$ over a field of characteristic zero is triangulable.

In this section we adopt the proof of this  result given in \cite[Ch. 5]{Es} to prove that every locally nilpotent $n$-graded derivation of $K[x,y]$ is triangulable by a tame $n$-graded automorphism. 

First recall some necessary definitions from \cite{Es}.

Let $0\neq w=(w_1,w_2)\in\mathbb{Z}^2$. Then $w$-degree of the monomial $x^{a_1}y^{a_2}$ is defined by $w(x^{a_1}y^{a_2})=a_1w_1+a_2w_2$. This degree function leads to the $w$-grading 
\bes
K[x,y]=\sum_d W_d
\ees
of $K[x,y]$, where $W_d$ is the span of all monomials of $w$-degree $d$.

Let $T=cx^{a_1}y^{a_2}\partial_i$ be a monomial derivation of $K[x,y]$, where $i=1,2$. Set $(s,t)=(a_1,a_2)-e_i$, where $e_i$ is the $i$-th vector of the standard basis of $K^2$. Then 
$$T(x^{m_1}y^{m_2})\in Kx^{m_1+s}y^{m_2+t}$$
 for all $m_1,m_2$. We call $(s,t)$ the \textit{strength} of $T$. 

Every derivation $D$ is a linear combination of monomial derivations. Set 
\bes
\mathrm{supp}\,D=\{(s,t)\in\mathbb{Z}^2\,|\, D \, \textit{contains a term of strength} \, (s,t)\}.
\ees

Let us denote by $D(s,t)$ the sum of all terms in $D$ of strength $(s,t)$ and set
\bes
D_p=\sum_{sw_1+tw_2=p} D(s,t).
\ees
Obviously, 
\bes
D=\sum_p D_p
\ees
and this decomposition is called the {\em w-homogeneous} decomposition of $D$. If $p$ is maximal with $D_p\neq 0$ then $p$ is called the {\em $w$-degree of $D$} and is denoted by $wdeg\,D$. When $w=(1,1)$ $p$ is called the {\em degree of D} and is denoted by $deg\,D$.

It is easy to check \cite{Es} that  $D_pW_d\subset W_{p+d}$ for all $p,d\in\mathbb{Z}$.

Consider the grading (\ref{g}) of $K[x,y]$. 
Set $K[y]_1=K[x,y]_1\cap K[y]$. Every triangular $n$-derivation of $K[x,y]$ can be written as $f\partial_x+\a\partial_y$ where $f\in K[y]_1$ and $\a\in K$.

\begin{pr}\label{p1}
Let $D$ be a locally nilpotent $n$-graded derivation of $K[x,y]$. Then there exists a tame $n$-graded automorphism $\alpha$ of $K[x,y]$ and $f(y)\in K[y]_1$ such that 
\bes
\alpha^{-1}D\alpha=f(y)\partial_{x}.
\ees
\end{pr}
\Proof
Let $D$ be a locally nilpotent $n$-graded derivation of $K[x,y]$. According to Corollary 5.1.16 in \cite[p. 91]{Es},  the following three cases are possible: 

$(i)$ $D=f(y)\partial_{x}$, for some $f(y)\in K[y]$; 

$(ii)$ $D=f(x)\partial_{y}$, for some $f(x)\in K[x]$; 

$(iii)$ there exist $s_0, t_0\geq 0$ such that $(s_0,-1)$ and $(-1,t_0)$ belong to $supp \,D$ and, furthermore, $\mathrm{supp}\, D$ is contained in the triangle with vertices $(s_0,-1)$, $(-1,-1)$, $(-1,t_0)$.

{\em Case $(i)$.} If $D=f(y)\partial_{x}$ with $f(y)\in K[y]_1$ then set $\alpha=\mathrm{id}$. Obviously, the identity automorphism is an $n$-graded automorphism.

{\em Case $(ii)$.} If $D=f(x)\partial_{y}$ with $f(x)\in K[x]_1$ then set $\alpha=(y,x)$. Obviously $\a$ is a $n$-graded automorphism of $K[x,y]$ and 
$\alpha^{-1}D\alpha=f(y)\partial_x$ with $f(y)\in K[y]_1$.

{\em Case $(iii)$.} Suppose that we have $s_0, t_0\geq 0$ such that $(s_0,-1),(-1,t_0)\in \mathrm{supp}\,D$. This implies that $D$ contains differential  monomials of the form $x^{s_0}\partial_y$ and $y^{t_0}\partial_x$ with nonzero coefficients. Hence $s_0=1+nk$, $t_0=1+nl$, $k,l\in\mathbb{Z}$ since $x^{s_0}, y^{t_0} \in K[x,y]_1$. 

Let $L$ be the line passing through the points $(1+nk,-1)$ and $(-1,1+nl)$. The defining equation of $L$ is 
$$(nl+2)x+(nk+2)y=n^2kl+nk+nl=nM.$$ 
Set $w=(nl+2,nk+2)$ and $p=n^2kl+nk+nl$. Obviously $\mathrm{wdeg}\,D=p$ and $D_p$ is the highest homogeneous part of $D$ with respect to the $w$-degree. It is well known that the highest homogeneous part of a locally nilpotent derivation is locally nilpotent (see, for example \cite[p. 90]{Es}). Consequently, $D_p$ is a locally nilpotent $n$-graded derivation. 

We can write $D_p=gD_1$, where $D_1=a\partial_x+b\partial_y$ with $gcd(a,b)=1$. By Corollary 1.3.34 in \cite[p. 29]{Es}, $D_1$ is locally nilpotent and $D_1(g)=0$. By Proposition 1.3.46 in \cite[p. 31]{Es}, $D_1$ has a slice in $K[x,y]$, i.e., there exists $s\in K[x,y]$ such that $D_1(s)=1$. This implies that $a(0,0)\neq 0$ or $b(0,0)\neq 0$. Assume that $a(0,0)\neq0$. This means that $D_1$ has a term of the form $c\partial_x$, where $c\in K^*$. Since $(1+nk,-1)\in \mathrm{supp}\,D_p$ and $D_p=gD_1$ it follows that $D_1$ also has a term of the form $dx^r\partial_y$ with $r\geq0$ and $d\in K^*$. 
Moreover, $g$ and $D_1$ are $w$-homogeneous since $D_p$ is $w$-homogeneous. 
Therefore $\mathrm{supp}\,D_1$ is on the line passing through the points $(-1,0)$ and $(r,-1)$. Notice that this line does not contain any other points with integer coordinates. Hence $D_1=c\partial_x+dx^r\partial_y$.  Since $D_p$ is an $n$-graded derivation it follows that $g\in K[x,y]_1$ and $n|r$. 

We have $g\in \mathrm{Ker}\,D_1=K[y-\frac{d}{(r+1)c}x^{r+1}]$ since $D_1(g)=0$. Consequently,  $g=a(y-\frac{d}{(r+1)c}x^{r+1})^N$ for some $a\in K^*$ and $N\in\mathbb{N}$ since $g$ is $w$-homogeneous. So
\bes
D_p=a(y-\frac{d}{(r+1)c}x^{r+1})^N(c\partial_x+dx^r\partial_y),
\ees
where $a,c,d\in K^{\ast}$, $r\geq0$, and $N\in\mathbb{N}$. Obviously, $t_0=N$ and $s_0=(r+1)N+r$. 

Let $\alpha$ be the automorphism given by 
\bes
\alpha(x)=x, \, \alpha(y)=y-\frac{d}{(r+1)c}x^{r+1}.
\ees 
This is an elementary $n$-graded automorphism since $n|r$. Direct calculations give that 
\bes
\alpha^{-1}D_1\alpha=c\partial_x
\ees
and 
\bes
\alpha^{-1}D_p\alpha=acy^{t_0}\partial_x. 
\ees
Since $\a$ is $w$-homogeneous, $\alpha^{-1}D_p\alpha$ is the highest $w$-homogeneous part of $\alpha^{-1}D\alpha$. Thus $\a$ turns all points of $\mathrm{supp}\,D_p$ to one point $(-1,t_0)$. Consequently, $s_0(\alpha^{-1}D\alpha)<s_0(D)$. Leading an induction on $s_0(D)+t_0(D)$ we can conclude that the statement of the proposition is true. 
$\Box$

\section{Automorphisms of $K[x^n, x^{n-1}y,\ldots,xy^{n-1}, y^n]$}

\hspace*{\parindent}

As we noticed above, every $n$-graded automorphism of $K[x,y]$ induces an automorphism of $K[x^n, x^{n-1}y,\ldots,xy^{n-1}, y^n]$. In this section we prove the reverse of this statement.

\begin{lm}\label{b4} Let $p\in K[x,y]$. If $p^n\in K[x,y]_0$ then $p\in K[x,y]_i$ for some $i\in  \mathbb{Z}_n=\mathbb{Z}/n\mathbb{Z}$.
\end{lm}
\Proof Consider the standard grading 
\bes
K[x,y]=A_0\oplus A_1\oplus\ldots\oplus A_k\oplus\ldots, 
\ees
where $A_i$ is the linear span of monomials of degree $i$ for all $i\geq 0$. For any $f\in K[x,y]$ denote by $f_i\in A_i$ its homogeneous part of degree $i$. 
Let 
\bes
p=p_{i_1}+p_{i_2}+\ldots+p_{i_k}, \ \ 0\neq p_{i_j}\in A_{i_j}, \ \ i_1<i_2<\ldots<i_k. 
\ees
Suppose that $p_{i_1},p_{i_2},\ldots,p_{i_s}\in K[x,y]_i$ for some $i\in  \mathbb{Z}_n$ and $p_{i_{s+1}}\notin K[x,y]_i$. Set $q=p_{i_1}+p_{i_2}+\ldots+p_{i_s}$. Obviously, $q^n\in K[x,y]_0$. Set $t=(n-1)i_1+i_{s+1}$. Then $t\not\equiv 0 \mod{n}$. We get  
\bes
(p^n)_t=(q^n)_t+np_{i_1}^{n-1}p_{i_{s+1}}=np_{i_1}^{n-1}p_{i_{s+1}}\notin  K[x,y]_0
\ees
since $q^n\in K[x,y]_0$. This contradicts to $p^n\in K[x,y]_0$. $\Box$

\begin{theor}\label{t1}
Every automorphism of $K[x^n, x^{n-1}y,\ldots,xy^{n-1}, y^n]$ over an algebraically closed field $K$ of characteristic zero is induced by an $n$-graded automorphism of $K[x,y]$. 
\end{theor}
\Proof Consider the derivation $D=y\partial_x$ of $K[x,y]$. Let $\overline{D}$ be the derivation of $K[x^n, x^{n-1}y,\ldots,xy^{n-1}, y^n]$ induced by $D$.

Let $\a$ be an arbitrary automorphism of $K[x^n, x^{n-1}y,\ldots,xy^{n-1}, y^n]$. Set $T=\a\overline{D}\a^{-1}$. This derivation is locally nilpotent since $D$ is locally nilpotent. Let $\widetilde{T}$ be the extension of $T$ to a derivation of $K[x,y]$ that uniquely defined by Lemma \ref{l1}. By Lemma \ref{l2}, $\widetilde{T}$ is a locally nilpotent $n$-graded derivation of $K[x,y]$. By Proposition \ref{p1}, there exists an $n$-graded tame automorphism $\b$ of $K[x,y]$ such that $S=\b^{-1}\widetilde{T}\b$ is a triangular $n$-graded derivation of $K[x,y]$. Let 
\bes
S=\b^{-1}\widetilde{T}\b=g(y)\partial_x, 
\ees
where $g(y)\in K[y]_1$. 
We get 
\bes
S(f)=g(y)\frac{\partial f}{\partial x}, \ \ f\in K[x,y]. 
\ees

Let $\overline{\b}$ be the automorphism of $K[x^n, x^{n-1}y,\ldots,xy^{n-1}, y^n]$ induced by $\b$. Then $S$ induces the derivation $\overline{S}=\overline{\b}^{-1}T\overline{\b}=\overline{\b}^{-1}\a\overline{D}\a^{-1}\overline{\b}$ of $K[x^n, x^{n-1}y,\ldots,xy^{n-1}, y^n]$. 

Let $\phi=\overline{\b}^{-1}\a$. Assume that $\phi(x^{n-i}y^i)=f_i$, where $0\leq i\leq n$. Applying the equation 
$\phi\overline{D}=\overline{S}\phi$ to $x^{n-i}y^i$ for all $i$, we get 
\bes
(n-i) f_{i+1}=g(y) \frac{\partial f_i}{\partial x}, 
\ees
i.e., 
\bes
0=g(y) \frac{\partial f_n}{\partial x}, f_n=g(y) \frac{\partial f_{n-1}}{\partial x},\ldots, (n-1)f_2=g(y) \frac{\partial f_1}{\partial x}, nf_1=g(y) \frac{\partial f_0}{\partial x}.
\ees
These equalities immediately give that 
\bes
\deg_xf_n=0, \deg_xf_{n-1}=1,\ldots,\deg_xf_{n-i}=i, \ldots, \deg_xf_0=n. 
\ees
In particular, $f_n\in K[y]$. 

We have 
\bee\label{f4}
\frac{f_0}{f_1}=\frac{f_1}{f_2}=\ldots=\frac{f_{n-1}}{f_n}
\eee
since the generators $x^n, x^{n-1}y,\ldots,xy^{n-1}, y^n$ of $K[x^n, x^{n-1}y,\ldots,xy^{n-1}, y^n]$ satisfy the relations 
\bes
\frac{x^n}{ x^{n-1}y}=\frac{ x^{n-1}y}{ x^{n-2}y^2}=\ldots=\frac{xy^{n-1}}{y^n}=\frac{x}{y}. 
\ees

Let $\frac{f_0}{f_1}=\frac{p}{q}$, where $p,q\in K[x,y]$ are relatively prime. Then $\frac{f_0}{f_n}=\frac{p^n}{q^n}$ by (\ref{f4}). Since $p^n$ and $q^n$ are relatively prime it follows that 
$f_0=p^nu$ and $f_n=q^nu$ for some $u\in K[x,y]$. Moreover, (\ref{f4}) implies that $f_i=p^{n-i}q^i u$ for all $i$. From this we get 
\bes
K[x^n, x^{n-1}y,\ldots,xy^{n-1}, y^n]\subseteq K+(u), 
\ees
where $(u)$ is the ideal of $K[x,y]$ generated by $u$. 
 This is possible if the leading word of $u$ divides all of the words $x^n, x^{n-1}y,\ldots,xy^{n-1}, y^n$. Consequently, $u\in K^*$. Over an algebraically closed field we can write $u=v^n$ for some $v\in K^*$. Replacing $vp$ with $p$ and $vq$ with $q$,  we may assume that $u=1$ and $f_i=p^{n-i}q^i$ for all $i$.

We have $q\in K[y]$ since $f_n=q^n\in K[y]$. We also have $\deg_x(p)=1$  since $p^n=f_0$  and $\deg_x(f_0)=n$. Set $p=xa(y)+b(y)$. We get 
\bes
K[x^n, x^{n-1}y,\ldots,xy^{n-1}, y^n]\subseteq K[f_n]+(p)\subseteq K[y]+(p), 
\ees
where $(p)$  is the ideal of $K[x,y]$ generated by $p$. 
Consequently, 
\bes
x^n=(xa(y)+b(y)) h +f(y). 
\ees
 Introducing a monomial order with prioritized $x$, we get that it is possible only if $a(y)=a\in K^*$. Consequently, $p=ax+b(y)$. By Lemma \ref{b4}, it implies that $p\in K[x,y]_1$ since $p^n\in K[x^n, x^{n-1}y,\ldots,xy^{n-1}, y^n]$. Set $\g=(ax+b(y),y)$. Then $\g$ is an elementary $n$-graded automorphism of $K[x,y]$. Set $\psi=\overline{\g}^{-1}\phi$. Then $\psi(x^{n-i}y^i)=x^{n-i}q^i$ for all $i$. We have 
\bes
K[x^n, x^{n-1}y,\ldots,xy^{n-1}, y^n]\subseteq K[q^n]+(x),  
\ees
where $(x)$ is the ideal of $K[x,y]$ generated by $x$. It is possible only if $q^n=cy^n$ for some $c\in K^*$. Consequently, $q=ey$ for some $e\in K^*$ since $K$ is algebraically closed. 

Let $\delta=(x,ey)$. Then $\overline{\delta}^{-1}\psi=\mathrm{id}$, i.e., $\overline{\delta}^{-1}\overline{\g}^{-1}\overline{\b}^{-1}\a=\mathrm{id}$. Consequently, 
$\a=\overline{\b} \overline{\g}\overline{\delta}=\overline{\b\g\delta}$ is induced by a tame $n$-graded automorphism of $K[x,y]$. $\Box$

\section{Amalgamated free product structure of $\mathrm{Aut}\,K[x^n, x^{n-1}y,\ldots,xy^{n-1}, y^n]$}

\hspace*{\parindent}

Let $G_n$ be the group of all $n$-graded automorphisms of the polynomial algebra $K[x,y]$. 
\begin{lm}\label{l3}
The subgroup $G_n$ of $\mathrm{Aut}\,K[x,y]$ is generated by all linear automorphisms and all automorphisms of the type $(x-\alpha y^m,y)$, where $m=1+ns$ is a positive integer and $\alpha\in K$.
\end{lm}
\Proof For any $f\in K[x,y]$ denote by $\bar{f}$ its highest homogeneous part with respect the standard degree function $\deg$. 
Let $\phi=(f,g)$ be a  $n$-graded automorphism of the algebra $K[x,y]$ and suppose that $\deg\,f=k$ and $\deg\,g=l$. If $k+l=2$ then $\phi$ is a linear automorphism. 

Suppose that $k+l\geq3$. It is well known that $k|l$ or $l|k$ (see, for example \cite{Cohn,Es}). Assume that $l|k$. In this case we have $\bar f=\alpha\bar g^m$ for some $\alpha\in K^*$ and $m\in \mathbb{N}$.  Since $\bar f, \bar g\in K[x,y]_1$ it follows that $m=1+ns$ for some $s\geq 0$. 
In fact, let $\deg(\bar f)=1+np$  and $\deg(\bar g)=1+nq$. Then 
\bes
1+np=m(1+nq). 
\ees
Consequently, $m-1=np-mnq=ns$. 

Therefore  $(x-\alpha y^m, y)$ is an elementary $n$-graded automorphism of $K[x,y]$. 
We have  
\[
(f,g)\circ (x-\alpha y^m, y)=(f-\alpha g^m, g)=(f',g),
\]
 where $\deg(f')<\deg(f)$. Leading an induction on $k+l$ we may assume that $(f',g)$ satisfies the statement of the lemma. Then $(f,g)$ also satisfies the statement of the lemma. $\Box$

\begin{co}\label{c0}
Every $n$-graded automorphism of $K[x,y]$ is $n$-graded tame. 
\end{co}

An automorphism $\phi\in\mathrm{Aut}\,K[x,y]$ is called {\em n-graded triangular}  if it has the form
\bes
\phi=(\alpha x+f(y),\beta y),
\ees
where $0\neq\alpha,\beta\in K$ and $f(y)\in K[y]_1$. 

Let $T_n$ be the group of all $n$-triangular automorphisms of the polynomial algebra $K[x,y]$. 

\begin{co}\label{c1}
$G_n=\mathrm{GL}_2\,(K)\ast_{B} T_n$, where $B=\mathrm{GL}_2\,(K)\cap T_n$.
\end{co}
\Proof Lemma \ref{l3} says that $G_n$ is generated by $\mathrm{GL}_2$ and $T_n$. Consider (\ref{f0}). We have $\mathrm{GL}_2\subseteq \mathrm{Aff}_2$, $T_n\subseteq \mathrm{Tr}_2\,(K)$, and $B\subseteq C$. This means that every decomposition of an element of $G_n$ 
in the form 
\bes
g_1g_2\ldots g_k, 
\ees
where $g_i\in \mathrm{GL}_2\cup T_n$ for all $i$ and $g_i$ and $g_{i+1}$ do not belong together to $\mathrm{GL}_2$ or $T_n$ for all $i<k$, determined by the amalgated free product structure (\ref{f0}). 
 Consequently, 
\bes 
G_n=\mathrm{GL}_2\,(K)\ast_{B} T_n \subseteq \mathrm{Aff}_2\,(K)\ast_{C}\mathrm{Tr}_2\,(K).  \ \ \ \Box
\ees

\begin{co}\label{c2} Let $E=\{\epsilon \mathrm{id} | \epsilon^n=1, \epsilon\in K\}$. Then 
\bes
\mathrm{Aut}\,K[x^n, x^{n-1}y,\ldots,xy^{n-1}, y^n]\cong G_n/E.
\ees
\end{co}
\Proof Consider the homomorphism 
\bee
\psi: G_n\to\mathrm{Aut}\,K[x^n, x^{n-1}y,\ldots,xy^{n-1}, y^n] 
\eee
defined by $\psi(\a)=\overline{\a}$, where $\overline{\a}$ is the automorphism of $K[x^n, x^{n-1}y,\ldots,xy^{n-1}, y^n]$ induced by the $n$-graded automorphism $\a$ of $K[x,y]$. 

By Theorem \ref{t1}, $\psi$ is an epimorphism. Let $\a \in \mathrm{Ker}\,\psi$. Then 
\bes
\a(x)^{n-i}\a(y)^i=x^{n-i}y^i
\ees
for all $0\leq i\leq n$. This implies that $\a(x)=\epsilon x, \a(y)= \epsilon y$ for some $n$th root of unity $\epsilon\in K$. Consequently, $\a\in E$. Obviously, $E\subseteq \mathrm{Ker}\,\psi$. $\Box$

Let 
\bes
\overline {\mathrm{GL}_2\,(K)}=\mathrm{GL}_2\,(K)/E, \overline{T_n}=T_n/E, \overline B=B/E. 
\ees

\begin{theor}\label{t2}
$\mathrm{Aut}\,K[x^n, x^{n-1}y,\ldots,xy^{n-1}, y^n]\cong \overline{\mathrm{GL}_2\,(K)} \ast_{\overline B} \overline{T_n}$.
\end{theor}
\Proof
By Corollaries \ref{c1} and \ref{c2}, the group $\mathrm{Aut}\,K[x^n, x^{n-1}y,\ldots,xy^{n-1}, y^n]$ is generated by $\overline {\mathrm{GL}_2\,(K)}$ and $\overline{T_n}$.

Let $G$ be any group and $\psi_1: \overline {\mathrm{GL}_2\,(K)} \to G$ and $\psi_2: \overline{T_n} \to G$ be any homomorphisms with $\psi_1|_{\overline{B}}=\psi_2|_{\overline{B}}$. 

Let $\a : {\mathrm{GL}_2\,(K)} \to\overline {\mathrm{GL}_2\,(K)}$ and $\b : T_n \to\overline{T_n}$ 
be natural homomorphisms. Set $\phi_1=\psi_1\a : {\mathrm{GL}_2\,(K)}\to G$ and $\phi_2=\psi_2\b : \mathrm{Tr}_n\to G$. Obviously, 
$\phi_1|_{B}=\phi_2|_{B}$. By the universal property of the amalgamated free products of groups \cite[Ch. 1]{Serre}, there exists a unique homomorphism $\phi: {\mathrm{GL}_2\,(K)}\ast_{B}T_n\to G$ such that $\phi|_{{\mathrm{GL}_2\,(K)}}=\phi_1$, $\phi|_{T_n}=\phi_2$.  Since $E\subseteq \mathrm{Ker}(\phi)$, $\phi$ induces the homomorphism $\overline{\phi} : ({\mathrm{GL}_2\,(K)}\ast_{B}T_n)/E\to G$. Obviously, $\overline{\phi}|_{\overline {\mathrm{GL}_2\,(K)}}=\psi_1$ and $\overline{\phi}|_{\overline {T_n}}=\psi_2$. By the definition of the amalgamated free product \cite{Serre}, we get 
\bes
({\mathrm{GL}_2\,(K)}\ast_{B}T_n)/E\cong \overline {\mathrm{GL}_2\,(K)} \ast_{\overline {B}} \overline{\mathrm{Tr}_n}. 
\ees
Corollary \ref{c1} finishes the proof of the theorem. 
 $\Box$\\

Recall that an automorphism $f\in\mathrm{Aut}\,K[x^n, x^{n-1}y,\ldots,xy^{n-1}, y^n]$ is called \textit{linearizable} if there exists $\phi\in\mathrm{Aut}\,K[x^n, x^{n-1}y,\ldots,xy^{n-1}, y^n]$ such that $\phi^{-1}f\phi$ is linear.

\begin{co}\label{c3}
Any automorphism of $K[x^n, x^{n-1}y,\ldots,xy^{n-1}, y^n]$ of finite order is linearizable.
\end{co}
\Proof By Corollary 1 in \cite[page 6]{Serre}  every element of
$\mathrm{Aut}\,K[x^n, x^{n-1}y,\ldots,xy^{n-1}, y^n]$ of finite order is conjugate to an element of $\overline {\mathrm{GL}_2\,(K)}$ or $\overline{T_n}$. Since $\overline{T_n}$ has no elements of finite order over a field of characteristic zero, any automorphism of $K[x^n, x^{n-1}y,\ldots,xy^{n-1}, y^n]$ of finite order is conjugate to an element of $\overline {\mathrm{GL}_2\,(K)}$. 
$\Box$

\bigskip

\begin{center}
	{\large Acknowledgments}
\end{center}

\hspace*{\parindent}

The second author is grateful to Max-Planck Institute f\"ur Mathematik for
hospitality and excellent working conditions, where part of this work has been done.

The second author is supported by the grant of the Ministry of Education and Science of the Republic of Kazakhstan (project  AP09261086).  

We are grateful to Professors M. Gizatullin, L. Makar-Limanov, and M. Zaidenberg  for their helpful suggestions and comments.

\end{document}